\documentclass{article} 
\usepackage{a4wide,epsfig,latexsym,xspace,oupbib,hyperref,url,xspace,amsmath}

\newcommand{\eg}{{\em e.g.\/}\xspace}

\newcommand{\etc}{{\em etc.\/}\xspace}

\newcommand{\HIDE}[1]{ }
\newcommand{\COMMENT}[1]{ }

\newcommand{\pc}{{\sc pc}\xspace}

\newcommand{\secref}[1]{\mbox{\S$\,$\ref{sec:#1}}}

\newcommand{\figref}[1]{\mbox{Figure~\ref{fig:#1}}}

\newcommand{\tabref}[1]{\mbox{Table~\ref{tab:#1}}}

\newcommand{\itref}[1]{\mbox{\ref{it:#1}}}

\newcommand{\cd}{\,|\,}

\newcommand{\E}{{\mbox{E}}}

\newcommand{\cip}{\mbox{$\perp\!\!\!\perp$}}

\newtheorem{expl}{Example}[section]
\newtheorem{definer}{Definition}[section]

\newtheorem{algor}{Algorithm}[section]

\newtheorem{rem*}{Remark}[section]

\newcommand{\halm}{\hspace*{\fill} $\Box$\par}

\newenvironment{ex}{\begin{expl}\rm}{\halm\end{expl}}

\newcommand{\pr}[1]{{p}(#1)}

\renewcommand{\pc}{\mbox{PC}\xspace}
\newcommand{\ec}{\mbox{EC}\xspace}
\renewcommand{\pr}{\mbox{\rm Pr}}
\newcommand{\bm}[1]{\mbox{$\mathbf #1$}}

\newcommand{\bR}{{\bm R}}

\newcommand{\rr}{\mbox{RR}\xspace}

\newcommand{\indo}[2]{\mbox{$#1 \,\cip\, #2$}}
\newcommand{\ind}[3]{\mbox{$#1 \, \cip\, #2 \cd
    #3$}}

\author{Philip~Dawid\thanks{University of Cambridge} \and
  Monica~Musio\thanks{University of Cagliari}\and
  Rossella~Murtas\thanks{University of Cagliari} }

\title{The Probability of Causation} 

\date{Dedicated to the memory of Stephen Elliott Fienberg\\[1ex]
  27 November 1942--14 December 2016\\[2.5ex]
  \today}

\begin{document}
\maketitle

\begin{abstract}
  \noindent Many legal cases require decisions about causality,
  responsibility or blame, and these may be based on statistical data.
  However, causal inferences from such data are beset by subtle
  conceptual and practical difficulties, and in general it is, at
  best, possible to identify the ``probability of causation'' as lying
  between certain empirically informed limits.  These limits can be
  refined and improved if we can obtain additional information, from
  statistical or scientific data, relating to the internal workings of
  the causal processes.  In this paper we review and extend recent
  work in this area, where additional information may be available on
  covariate and/or mediating variables.\\

  
\end{abstract}

\section{Introduction}
\label{sec:intro}
Many legal proceedings hinge on the assignment of blame or
responsibility for some undesirable outcome. In a civil case a patient
may sue a pharmaceutical company for damage caused as a side effect of
one of its products; or a state health department may bring an action
against a tobacco company for not disclosing information it had on the
health risks of smoking, thus (it is claimed) leading to unnecessary
deaths.

In many such cases there will be no dispute about the facts at issue.
The patient took the drug and suffered the side effect.  The tobacco
company admits to having withheld the information.  What remains at
issue is the causal relation between the established facts.  But any
attempt to understand this quickly throws us into the philosophical
quagmire of counterfactual reasoning, where we have to consider what
might have happened in circumstances known to be false.  Whether or
not the patient's suit succeeds will depend on whether or not she can
prove, to the appropriate legal standard (\eg, ``on the balance of
probabilities'') that she would not have developed the same outcome,
had she not taken the pharmaceutical product.  The damages that the
tobacco company are liable for will depend on an assessment of how
many lives could have been saved, had they made the information
available.

In this article we concentrate on cases of the first kind, where it is
desired to assess whether or not the same outcome would have occurred
had the putative causal event been absent.  Typically, epidemiological
evidence will be admitted as to the frequency of the adverse outcome,
both in patients who have, and in patients who have not, been exposed
to the product.  Although such evidence is clearly of relevance, it is
less clear exactly how.  A common procedure is to compute the {\em
  relative risk\/} (\rr), obtained by dividing the frequency of the
outcome among those exposed to the product by the corresponding
frequency among those not so exposed.  And it is frequently asserted
that a relative risk exceeding 2 is enough to prove that the
``probability of causation'' (\pc) exceeds 0.5, and thus to establish
a causal link ``on the balance of probabilities''.

In this paper we consider such arguments in greater detail.  We
explain why it is difficult to establish a precise value for the
probability of causation on the basis of epidemiological or other
scientific evidence, which at best can only provide interval bounds
for \pc.  We also show how it is possible to reduce the uncertainty
about \pc by collecting and taking account of additional data, so
shedding some light on the internal workings of the causal black box.

The work presented here builds on a long-standing interaction and
collaboration with Stephen Fienberg, as represented in particular by
the recent papers
\textcite{sef/dlf/apd:socmeth,sef/dlf/apd:pearlreply,apd/mm/sef:ba}.
This formed a very small component part of Steve's many highly
significant contributions to the correct use of Statistics in the Law,
which work, substantial as it was, itself constituted only a very
small component of his numerous original and highly influential
contributions over an enormous range of statistical topics.  He was a
delightful friend and a stimulating colleague and companion; he is
sorely missed.

\section{Group to individual analysis: The basic problem}
\label{sec:g2i}
\begin{figure}[htbp]
  \centering
  \includegraphics[width=.5\linewidth]
  {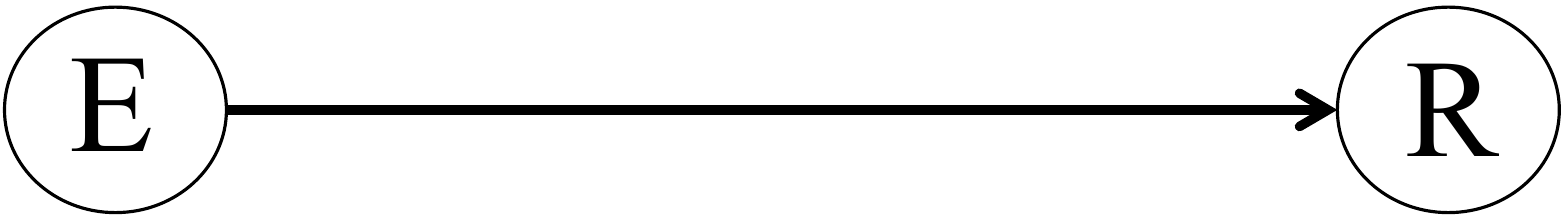}
  \caption{Simple cause-effect relationship}
  \label{fig:basic}
\end{figure}

Epidemiology attempts to discover causal relationships, such as
between a certain exposure, $E$ and a certain response, $R$.  Such a
relationship can be intuitively represented by a diagram such as that
of \figref{basic}.  (There is a formal semantics underlying such
diagrams, but we shall not delve into that here).

Epidemiological experts are frequently called on to testify in cases
revolving around the assignment of individual responsibility.
However, such testimony can only concern the behaviour of groups of
individuals, and its relevance to any individual case is at best
indirect.  The general issue of ``group-to-individual'' (``G2i'')
inference has recently become a centre of attention both within and
outside the Law \cite{g2i}.

In our context, the important distinction to be made is between
inference about ``Effects of Causes'' (``EoC''), and inference about
``Causes of Effects'' (``CoE'').  The former, concerning the likely
outcome of a contemplated intervention or exposure, in a specific or
generic individual, is the focus of most scientific investigation and
experimentation, including epidemiological studies.  However, it is
the latter, concerning the possible causes of an observed effect in an
individual case, that is of most relevance for the Law.  The
distinction between EoC and CoE is crucial, but all too often goes
unremarked and unappreciated.

We illustrate the subtle relationship between these two forms of
inference by means of an example
\cite{apd:aberdeen,sef/dlf/apd:socmeth}.

\begin{description}
\item[Effects of Causes (EoC)] Ann has a headache.  She is wondering
  whether to take aspirin.  Will that cause her headache to disappear?
\item[Causes of Effects (CoE)] Ann had a headache and took aspirin.
  Her headache went away.  Was that caused by the aspirin?
\end{description}

To assist in addressing these questions, we assume that we have data
from a very large, well-conducted, randomised trial of the success of
aspirin in curing headaches, whose subjects can be assumed to be
similar to Ann in the way they respond.  This yields percentages as in
the following table:
\begin{table}[htbp]
  \centering
  \begin{tabular}[t]{lccc}
    & Recover &Not recover & Total\\
    No aspirin & 12 & 88 & 100\\
    Aspirin & 30 & 70 & 100
  \end{tabular}
  \caption{Results of aspirin trial}
  \label{tab:aspirin}
\end{table}

\section{A framework for analysis}
\label{sec:framework}

\begin{description}
\item[Effects of Causes (EoC)] Ann can argue as follows ``If I take
  the aspirin, I will be just like the subjects in the treatment
  group, so I can assess my probability of recovery in this case as
  30\%.  Similarly, If I don't take the aspirin, my probability of
  recovery in that case can be assessed to be 12\%.  Other things
  being equal, it would be better to take the aspirin.''  While this
  argument does, not, perhaps address the question exactly as posed,
  it does give Ann fully adequate guidance to address and solve her
  EoC decision problem.  What is at issue was the as yet unknown {\em
    event\/} of recovery, and the data are directly relevant for
  assessing uncertainty about that unknown outcome (for each of the
  two decisions that Ann might take).

\item[Causes of Effects (CoE)] Things are very different for the CoE
  question.  Indeed, there is now no unknown {\em event\/}: We know
  that Ann definitely took the aspirin, and she definitely recovered.
  What is at issue here is not uncertainty about an as yet unknown
  event, but uncertainty about a {\em relationship\/} between known
  events: the existence, or otherwise, of a causal relationship
  between taking the aspirin and recovering.
\end{description}

\subsection{Counterfactuals and potential outcomes}
\label{sec:po}
One way of thinking about whether a certain outcome was {\em caused\/}
by a certain exposure is to ask the question: ``Would the same outcome
still have occurred, even if the exposure had been absent?''.  Because
this question is posited on a hypothesis that is known to be false,
and thus counter to known facts, it is termed a counterfactual query.
There is much discussion of counterfactuals in the philosophical
literature (see \eg\ \textcite{menzies:sep}), and they raise many
perplexing problems.  In particular, is it even meaningful to regard a
counterfactual as having a definite (even if unknown) truth value?

Supposing that we can talk of the truth value of counterfactual, we
might reinterpret the assertion ``The exposure caused the outcome''
as:

\begin{quote}
  ``Both the exposure and the outcome occurred; and had
  (counterfactually) the exposure not occurred, neither would the
  outcome have occurred.''
\end{quote}

With this interpretation, and taking account of the known facts, the
truth value of the assertion ``The exposure caused the outcome''
becomes identified with the (typically unknown) truth value of the
counterfactual proposition ``In the counterfactual world where the
exposure did not occur, the outcome did not occur''.

\paragraph{Some comments}
\begin{enumerate}
\item There is an implicit ``ceteris paribus'' condition in the above
  interpretation.  In the counterfactual world being considered, we
  are changing the truth value of the exposure, but nothing else---or,
  at any rate, nothing else occurring before, or simultaneously with,
  the exposure.  We must of course allow for the outcome to change,
  and we might possibly also allow changes to other, later, events
  (perhaps intermediate between the exposure and the outcome), in the
  light of the change to the exposure.  Since there are various ways
  in which our counterfactual comparison world can be constructed,
  there is a resultant ambiguity in the truth value of the
  counterfactual assertion.  \textcite{dkl:book} argues that we should
  use the ``closest possible world'' to the real world, subject to the
  counterfactual change in the exposure---but this still leaves a
  great element of vagueness.
\item There may be some choice as to what the counterfactual value of
  the exposure would have been, if it had been other than it actually
  was.  If Ann had not taken two aspirins, she might have taken none
  at all, or taken three aspirins, or taken two paracetamol tablets,
  or gone to bed with a wet towel around her head.  The counterfactual
  query only becomes precise when we have fully specified the
  ``counterfactual foil'' for comparison with the factual exposure.
\item There will typically be a variety of different exposures that we
  can consider.  For example, Ted has been a smoker for 20 years, has
  lived by a busy main road which is subject to high air pollution,
  and has developed lung cancer.  Then we could ask ``Was Ted's cancer
  caused by his smoking?'', ``Was it caused by his air pollution
  exposure?'', ``Was it caused by his smoking and the air pollution
  jointly?'', \etc, \etc\@ Each of these questions involves
  consideration of a different counterfactual world for comparison
  with the real world, and there is no reason why we could not regard
  all of them as true simultaneously---which leads to some
  difficulties in understanding what might be the ``actual cause'' of
  Ted's cancer (see \textcite{halpern} for a thoughtful analysis of
  the elusive concept of ``actual cause'').  In particular, we should
  not express ``Ted's smoking caused his cancer'' as ``The cause of
  Ted's cancer was his smoking''---though we might say ``A cause of
  Ted's cancer was his smoking''.  The former phrasing suggests,
  misleadingly, that there is a variable, ``The cause of Ted's
  cancer'', that can take different values (\eg, smoking, air
  pollution, natural causes).  But as these putatuve causes are not
  mutually exclusive, there can be no such variable.  And when we can
  assign probabilities to ``Ted's smoking caused his cancer'', ``Air
  pollution caused Ted's cancer'', \etc, it is perfectly acceptable
  for these to sum to more than 1.  (For a specific example of this,
  see \textcite{hutchinson:etal:91b}).
\end{enumerate}

The standard framework of probability theory does not easily admit
argument about counterfactual worlds and uncertain relationships.  To
formulate such questions meaningfully, we need to expand the universe
of discourse.  The way this has traditionally been done in statistics
is by the introduction of so-called ``potential outcomes''
\cite{jn:jrss,dbr:as}.  We consider a binary exposure, or cause,
variable $E$, and a binary response, outcome, or effect, variable $R$.
For each value $e = 0$ or $1$ of $E$, we introduce a variable $R_e$,
conceived as the value that $R$ would assume, were it the case that
$E= e$.  If in the actual world $E =1$ (say), then the observed
outcome will be $R=R_1$; in this case we have no way of observing
$R_0$, which thus becomes a (necessarily unobservable) counterfactual
variable. Similarly, if in the actual world $E =0$, then $R=R_0$, and
$R_1$ becomes counterfactual and unobservable. It is typically
assumed, explicitly or implicitly, that the bivariate variable
$\bR :=(R_0,R_1)$ is determined before the value of $E$ is, and does
not depend on how that value of $E$ came to be determined or known
(for example, whether by external intervention, or by observation
alone).  The effect of setting or observing $E=1$ (say) is just to
uncover the pre-existing value of $R_1$, and identify $R$ with
$R_1$. In particular, although the pair $(R_0,R_1)$ is modelled
probabilistically just like any other bivariate random variable, there
is the epistemic problem that there are no circumstances in which it
would be possible to observe both its constituent values.  This has
been termed ``the fundamental problem of statistical causality''
\cite{pwh:jasa}.

The framing and analysis of causal questions in terms of potential
outcomes has become the industry standard in statistical causality.
For the investigation of ``effects of causes’'', it has been strongly
argued that this approach is both unnecessary and potentially
misleading \cite{apd:cinfer,apd:annrev}.  However, for studying causes
of effects there is currently no alternative to the use of potential
outcomes to model counterfactual possibilities.  Nevertheless, the
logical difficulty that the bivariate quantity $\bR$ can never be
fully observed creates subtle problems that must not be ignored.

\subsection{The probability of causation}
\label{sec:pc}

Armed with the above machinery and notation, we can now formulate an
expression for the probability of causation.

We interpret the expression ``the exposure $E=1$ causes the outcome
$R=1$'' as the event ``\ec: $R_1 = 1, R_0 = 0$''.  This encodes the
property that the potential response to $E=1$ is $1$, while at the
same time the potential response to $E=0$ is $0$.  However, although
we have thus expressed the causality relation as an event, in what
appears to be a standard probabilistic framework, the fundamental
problem of statistical causality implies that this is an event whose
truth value can never be determined.

Suppose now that (for Ann) we have observed $E=1$ and $R=1$ (or,
equivalently, $E=1$ and $R_1=1$).  Since we must condition on all
known facts, the appropriate expression for the probability of
causation (in Ann's case) is
\begin{eqnarray}
  \nonumber
  \pc &=& \pr(\ec \mid E=1, R=1)\\
  \label{eq:pc}
      &=& \pr(R_0 = 0 \mid E=1, R_1=1). 
\end{eqnarray}
However, merely establishing a notational expression for the
probability of causation does not, of itself, solve the problem of how
to evaluate this quantity, involving as it does the probability of an
unobservable event.

The approach we shall take is to regard the pair $(R_0,R_1)$ as having
a joint distribution (jointly with all other variables in the
problem), but to make no assumptions about that distribution, other
than any constraints imposed by empirically observable information.
Typically it will not then be possible to identify that joint
distribution fully, but only to constrain it within limits---though
those limits will depend on the precise form of the available
information.

There is a further point that needs to be taken into account.  The
variables in \eqref{eq:pc} are those pertaining specifically to Ann,
and the joint distribution $\pr$ appearing in it is thus a
distribution over those specific variables.  But who is the analyst
assigning this distribution?  It might be Ann herself, or some
onlooker.  These individuals will typically have different background
information, and so different distributions.  In particular, an
onlooker might reasonably infer, merely from observing that Ann
decided to take the aspirin ($E=1)$, that her headache is particularly
bad, and consequently that her potential responses $R_0$ and $R_1$ are
both likely to be poorer than if he did not have this information.
For Ann, however, the severity of her headache is known background
information that is already taken into account in her own
distribution, and the mere fact that she is contemplating taking the
aspirin gives her no additional information about her potential
responses.  That is to say, for Ann's own probability distribution we
can reasonably invoke the ``no confounding'' assumption,
\begin{equation}
  \label{eq:noconfann}
  \indo E \bR,
\end{equation}
which expresses probabilistic independence between the fact of her
exposure, $E$. and the pair $\bR$ of her potential responses.  In this
case expression \eqref{eq:pc} simplifies to
\begin{equation}
  \label{eq:pc0}
  \pc = \pr(R_0 = 0 \mid  R_1=1). 
\end{equation}
However, for the onlooker who is not initially aware of the severity
of Ann's headache, it may be inappropriate to accept
\eqref{eq:noconfann}, and hence \eqref{eq:pc0}.

In most of the sequel (with the exception of \secref{allow}) we shall
assume that the analyst has sufficient background information about
Ann for \eqref{eq:noconfann} to be true (to a good enough
approximation), and hence use \eqref{eq:pc0} as the expression for the
probability of causation.  Furthermore, because we shall wish to use
external data from other individuals to assist in estimating the
analyst's probability distribution $\pr$ for Ann, we shall also need
to assume that the background information the analyst has on those
individuals is essentially the same as what he\footnote{We refer to
  the analyst in the masculine, to distinguish him from the subject,
  Ann---while not precluding that they could be the same individual.}
has for Ann --- a condition that can sometimes be met by restricting
the set of external individuals considered.  Even so, these
assumptions, necessary for much of our analysis, are highly
restrictive, and will often be inappropriate.  Careful consideration
must be given to their suitability before applying our results.

\section{Basic analysis}
\label{sec:basic}

We start by considering the case that the only available empirical
information is that presented in \tabref{aspirin}.

Assuming the analyst can regard the individuals in the trial as
similar to Ann, he might argue: ``When Ann takes aspirin, she puts
herself in the same position as those external individuals in the
treatment group, of which $\pr(R=1 \mid E=1) = 30\%$ recovered.
Because of the no confounding assumptions (both in the data and for
Ann), this means that her potential response to taking aspirin, $R_1$,
has a 30\% probability of taking the value $1$ (Recover).  I can thus
identify my marginal probability (for Ann), $\pr(R_1=1)= 0.3$.
Similarly, I can evaluate my marginal probability $\pr(R_0=1)=0.12$.
However, because of the fundamental problem of statistical causality,
there is no further empirical information available relating to my
joint distribution of $(R_0,R_1)$.''

Armed with this limited information, the analyst can attempt to fill
in the entries in \tabref{rr}, specifying his bivariate probability
distribution for $(R_0,R_1)$.  The margins are fully determined by the
above empirically informed probabilities, but since he has no
empirical information directly relevant to the joint probability
$\pr(R_0=0,R_1=1)$, he simply enters an unknown (and unknowable) value
$x$ in the associated internal cell of the table.  In combination with
the known marginal probabilities, this allows completion of the body
of the table, as shown.

\begin{table}[htbp]
  \centering
  \caption{Joint distribution of potential responses $R_0$ and $R_1$.
  }
  \begin{tabular}[t]{c|cc|c}
    \multicolumn{1}{c}{} &\multicolumn{2}{c}{$R_0$}\\
    $R_1$ & 0 & 1 \\
    \hline
    0 &  $0.88-x$ & $x-0.18$  & 0.70\\
    1 & $x$ & $0.30-x$  & 0.30\\
    \hline
                         & 0.88 & 0.12 & 1
  \end{tabular}
  \label{tab:rr}
\end{table}

Now although the precise value of $x$ remains unknown, from the table
we can extract partial information about it, using the fact that any
probability must be nonnegative.  Applying this principle to the four
entries in the body of the table, we obtain: $x \leq 0.88$,
$x \geq 0.18$, $x\geq 0$, $x \leq 0.30$.  Combining these and omitting
redundancies gives the inequality: $0.18 \leq x \leq 0.30$.  This
interval bound is the full information available about $x$ on the
basis of the empirical data.

Now because of assumption \eqref{eq:noconfann},
$\pc = \pr(R_0=0 \mid R_1 = 1) = x/0.30$.  Thus the analyst can
assert:
\begin{equation}
  \label{eq:x30}
  0.6 \leq  \pc \leq 1.
\end{equation}
Even under the restrictive assumptions that have been made, this wide
interval bound is all that can be deduced, from the empirical data,
about the probability of causation in this particular case.

The above logic can be applied in the general case, leading to the
following interval bound in terms of empirically estimable quantities
(under our no-confounding assumptions):
\begin{equation}
  \label{eq:rrA}
  \max\left\{0,1 - \frac{1}{\rr}\right\}
  \leq
  \pc 
  \leq
  \min\left\{1,  \frac{\pr(R=0 \mid E=0)}{\pr(R=1 \mid E=1)}\right\},
\end{equation}
where
\begin{equation}
  \label{eq:exptrre}
  \rr := \frac{\pr(R = 1 \mid E=1)}{\pr(R = 1 \mid E=0)}
\end{equation}
is the {\em (causal) risk ratio\/}.  We note that $\rr > 2$ implies
$\pc > 0.5$, which might be taken to prove causation ``on the balance
of prbabilities''.  However, if $\rr < 2$ we can {\em not\/} deduce
from \eqref{eq:rrA} that $\pc < 0.5$ (to reject causation, on the
balance of probabilities).

\section{Allowing for unobserved confounding}
\label{sec:allow}
When we cannot assume \eqref{eq:noconfann} we cannot simply apply the
above theory directly, since the very fact that Ann decided to take
the aspirin gives the analyst some indirect information about Ann's
state of health---in the light of which it would no longer be
appropriate to consider her as similar to the individuals in the
experimental study, for which this information was not available.
\textcite{tian/pearl:probcaus} showed that in this situation we can
still obtain bounds for $\pc$, so long as, in addition to the
experimental data, we also have observational data on individuals
whose behaviour (in the opinion of the analyst) can be considered
similar to that of Ann, in that they have the same dependence as she
does between their decision $E$ to take the aspirin and their pair
$\bR$ of potential responses.  \textcite{tian/pearl:probcaus} then
derive the following interval bounds:
\begin{equation}
  \max\left\{0, \frac{\pr(R=1) - \pr(R_0=1)}{\pr(E=1, R=1)}\right\}
  \leq \pc \leq
  \min\left\{1, \frac{1 -\pr(R_0 =1)-\pr(E=0,R=0)}{\pr(E=1,R=1)}\right\},
  \label{eq:tpint}
\end{equation}
where $\pr(R_0=1)$ denotes the recovery probability in the
experimental control group, who were externally assigned to have $E$
set to $0$; while the other probabilities are obtained from the
observational data (interestingly, we do not require experimental data
on those who were assigned treatment).

For example, suppose that, in addition to the experimental data of
\tabref{aspirin}, we have observational data in which
$\pr(E=0, R=1) = 0.12$.  We then find that the lower bound in
\eqref{eq:tpint} is $1$---allowing us to establish causation with
certainty.

\section{Peeking into the box}
\label{sec:peeking}
So far we have considered the relationship between exposure and
response as a causal ``black box'', without attempting to understand
or model its hidden mechanisms in any further detail.  If we are
willing to make some assumptions about these mechanisms, or access
external information about them, we can make better causal inferences
from our data---even in the absence of any further information about
Ann.

\subsection{Monotonicity and beyond}
\label{sec:mon}
One assumption commonly made is {\em monotonicity\/}.  In a context
such as ours, this would say that, if recovery would have occurred
without the treatment, it would certainly occur with the treatment.
That is to say, $R_0 = 1\Rightarrow R_1 =1$.  This implies in
particular $\pr(R_0=1)\leq \pr(R_1=1)$, or equivalently $\rr \geq 1$.
This is a constraint on the margins of a table such as
\tabref{aspirin}, and so could be falsified; but even when this
marginal constraint is satisfied, there is typically no further
information to help us decide whether monotonicity holds or not.
Nevertheless it is frequently taken as a reasonable requirement, and
would add one further constraint on the entries in \tabref{aspirin}:
$\pr(R_0 = 1, R_1=0) = 0$.  Since in that table
$\pr(R_0 = 1, R_1=0) = x-0.18$, we can then compute the exact value
$x = 0.18$, so that $\pc = 0.18/0.30 = 0.6$.  Here (as in general),
monotonicity allows as to replace the interval given by \eqref{eq:rrA}
by its lower endpoint $1 - 1/\rr$.  In this case, $\rr <2$ does imply
$\pc < 0.5$.

Monitonicity, while appealing, is a strong, perhaps overstrong,
assumption about the internal workings of the causal black box.
Moreover, because it is a property of the pair of values $(R_0,R_1)$
considered jointly, the fundamental problem of statistical causality
implies that we can never expect to obtain empirical evidence to let
us decide whether or not monitonicity does in fact hold.  For that
reason we take a different approach in the rest of this paper,
focusing on using additional, empirically observable, information to
shed some light on internal causal mechanisms, and thus improve causal
inferences.  In such cases, although we typically cannot obtain an
exact point value for the probability of causation, we may be able to
narrow the interval bounds on it.  We shall particularly consider
cases in which we might have information on additional covariates
and/or mediators that influence the causal pathway.

\section{Covariate information}
\label{sec:cov}

\begin{figure}[htbp]
  \centering
  \includegraphics[width=.5\linewidth]
  {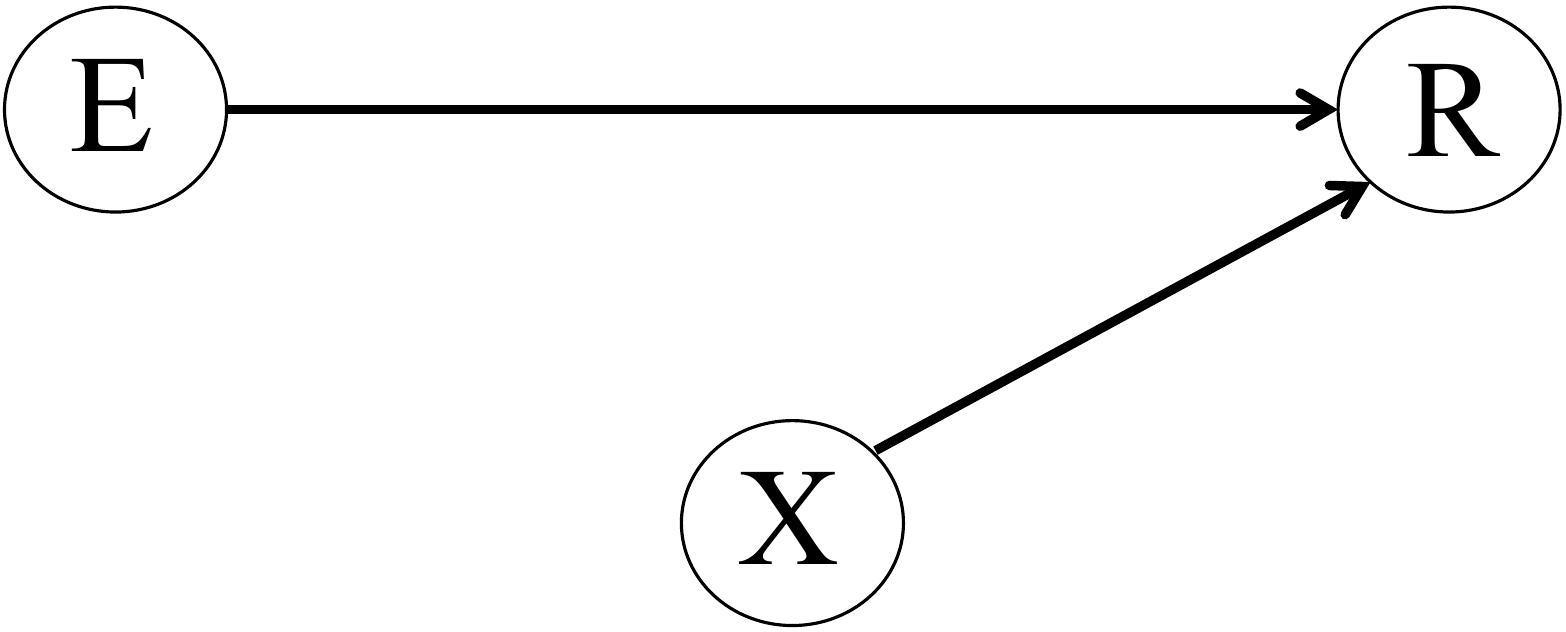}
  \caption{Covariate}
  \label{fig:expcov}
\end{figure}

Suppose that, in the experimental data, we can also observe an
additional covariate, $X$ --- an individual characteristic that can
vary from person to person, and can affect that person's response.
This situation is represented by the diagram of \figref{expcov}.
For simplicity, we suppose that $X$ can only take a discrete set of
values.  We can then estimate, from the data, the response
probability, $\pr(R=1 \mid E=e, X=x)$, among those individuals having
covariate value $X=x$ and exposure level $e$ ($=0$ or $1$).

In those cases that we are also able to measure $X$ for Ann, say
$X=x^*$, we can simply restrict the experimental subjects to those
having the same covariate value (who are thus like Ann in all relevant
respects), and apply the analysis of \secref{basic}.  We can then use
formula \eqref{eq:rrA}, after first replacing all the probabilities in
\eqref{eq:rrA} and \eqref{eq:exptrre} by their conditional versions,
given $X=x^*$.

However, it turns out that we can make use of the additional covariate
information in the experimental data, even when we cannot observe $X$
for Ann.  For this case, it is shown in \textcite{apd:aberdeen} that
we have the interval bounds
\begin{equation}
  \label{eq:better}
  \frac{\Delta}{\Pr(R = 1\mid E=1)} \leq \pc \leq 1 -
  \frac{\Gamma}{\Pr(R = 1\mid E=1)}
\end{equation}
where
\begin{eqnarray*}
  \Delta &=& \sum_x  \Pr(X=x) \times\max\left\{0, \Pr(R= 1
             \mid E=1, X=x) - \Pr(R = 1 \mid E=0, X=x)\right\}\\
  \Gamma &=&  \sum_x \Pr(X=x) \times \max\left\{0, \Pr(R = 1 \mid E=1,
             X=x) - \Pr(R = 0 \mid E=0, X=x)\right\},
\end{eqnarray*}
and that this interval is always contained in that given by
\eqref{eq:rrA}, which ignores the additional covariate information.

\begin{ex}
  Suppose that the covariate $X$ is binary, taking values $0$ and $1$
  with equal probability.  Suppose further that, from the experimental
  data, we obtain the following probabilities:
  \begin{eqnarray*}
    \pr(R=1 \mid E=1, X=0) &=& 0.12\\
    \pr(R=1 \mid E=0, X=0) &=& 0.24\\
    \pr(R=1 \mid E=1, X=1) &=& 0.60\\
    \pr(R=1 \mid E=0, X=1) &=& 0.12.
  \end{eqnarray*}
  These imply marginal recovery probabilities (not taking the value of
  $X$ into account):
  \begin{eqnarray*}
    \pr(R = 1\mid E=1) &=& 0.36\\ 
    \pr(R = 1\mid E=0) &=& 0.18.
  \end{eqnarray*}
  \begin{enumerate}
  \item \label{it:nox} If we had not observed $X$, either in the data
    or in Ann, we would use the marginal values, to obtain a relative
    risk of $2$, and interval bounds
  $$0.5 \leq \pc \leq 1.$$
\item If we observe $X$ both in the data and in Ann, we obtain
  $$0.8 \leq \pc \leq 1$$
  if Ann has $X=1$, and
  $$0 \leq \pc \leq 1$$
  if Ann has $X=0$.
\item Finally, if we have observed $X$ in the data but not in Ann, the
  relevant interval becomes
  $$0.67 \leq \pc \leq 1,$$
  an improvement on that of \itref{nox}
\end{enumerate}
\end{ex}

\section{Sufficient covariate}
\label{sec:obsconf}
\begin{figure}[htbp]
  \centering
  \includegraphics[width=.5\linewidth]
  {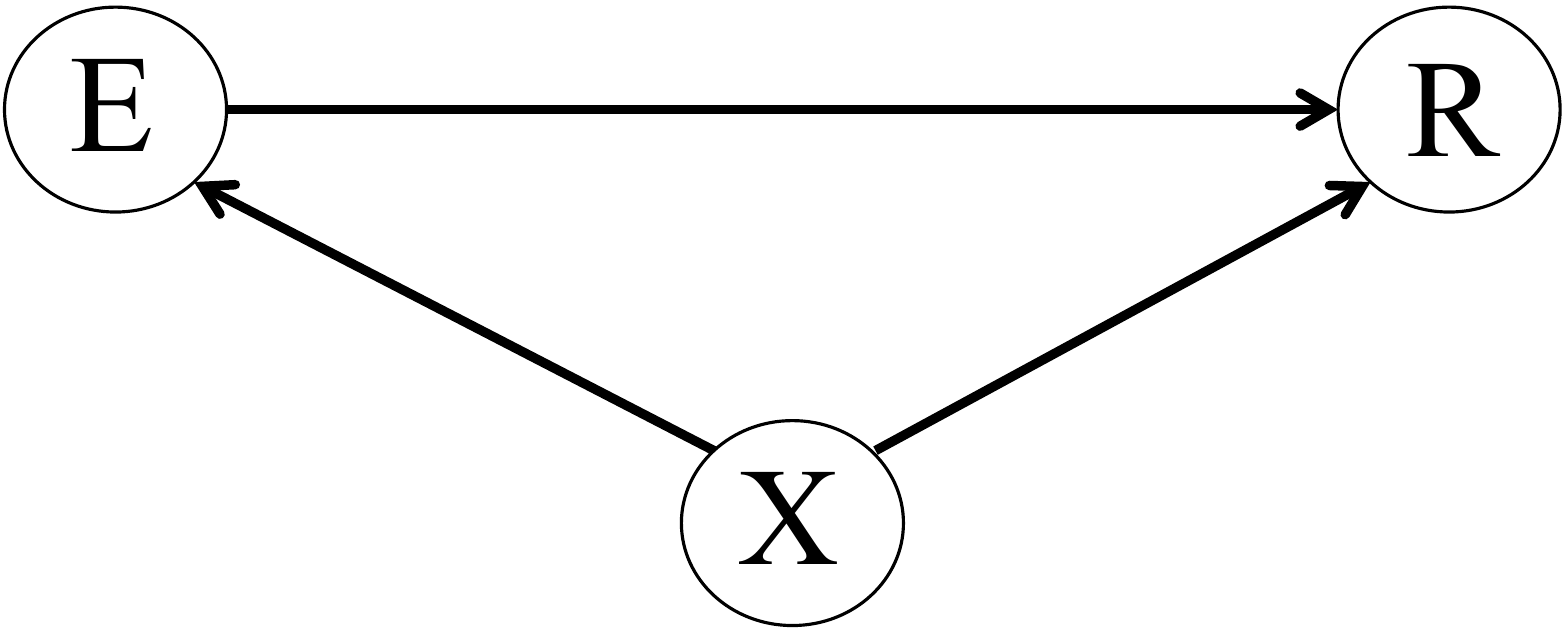}
  \caption{Sufficient covariate}
  \label{fig:suffcov}
\end{figure}

In \secref{cov} above it was assumed that the assignment of exposure
$E$ was unrelated to the covariate $X$.  Here we relax this, and allow
$X$ to influence both the exposure $E$ and the response $R$.  But we
assume that there are no further, unobserved, variables that might act
to confound this relationship---$X$ is a ``sufficient covariate''
\cite{hg/apd:aistats2010,apd:annrev} for assessing the causal effect
of $E$ on $R$.  The relevant diagram is now that of \figref{suffcov}.
We assume that the identical same structure holds for the individuals
in the data and for Ann.  $X$ is observed in the data, but possibly
not for Ann.

When $X$ is observed for Ann, say $X=x^*$, we can again apply the
basic inequalities \eqref{eq:rrA}, after replacing the probabilities
in \eqref{eq:rrA} and \eqref{eq:exptrre} by their conditional
versions, given $X=x^*$.

When $X$ is not observed for Ann, the interval bounds for $\pc$
become:
\begin{align}
  \nonumber
  \frac 1 {\pr(R=1 \mid E=1)}
          \times
          \sum_x & \max\{0, \{\pr(R=1 \mid E=1, X=x) - \pr(R=1 \mid E=0, X=x)\}\\
  \nonumber
       &{}\times \pr(X=x \mid E=1)\\[2ex]
  \nonumber
       &\leq \pc \leq\\[1ex]
  \nonumber
  1-\frac 1 {\pr(R=1 \mid E=1)}
          \times\sum_x &\max\{0, \{\pr(R=1 \mid E=1, X=x) - \pr(R=0 \mid E=0, X=x)\}\\
       &{}\times \pr(X=x \mid E=1).
         \label{eq:suffcovint}
\end{align}

Now in this case we can alternatively use the ``back-door formula''
\cite{pearl:book}:
\begin{equation}
  \label{eq:backdoor}
  \pr(R_0=1) = \sum_x\pr(R=1 \mid E=0, X=x)\times\pr(X=x)
\end{equation}
to recreate the experimental response rate among the controls---and so
compute the interval \eqref{eq:tpint}.  However, this would be
tantamount to ignoring the additional information on $X$ in the data.
It can be shown that the interval given by \eqref{eq:suffcovint} is
always contained in that given by \eqref{eq:tpint}.

\begin{ex}
  \label{ex:suffcov}
  Suppose $X$ is binary, and from the data we obtain the following
  probabilities:
  \begin{eqnarray*}
    \Pr(X=1) &=& 0.5\\
    \Pr(E=1\mid X=0) &=& 0.8\\
    \Pr(E=1\mid X=1) &=& 0.2\\
    \Pr(R=1 \mid E = 1, X=1) &=& 0.2\\
    \Pr(R=1 \mid E = 0, X=1) &=& 0.8\\
    \Pr(R=1 \mid E = 1, X=0) &=& 0.8\\
    \Pr(R=1 \mid E = 0, X=0) &=& 0.2.
  \end{eqnarray*}
  Then we obtain the following lower bounds for the probability of
  causation (the upper bound being 1 in all cases):
  \begin{tabbing}
    When Ann is observed to have
    $X=1$:\quad\quad\quad\quad\quad\quad\quad\quad\quad\quad\= $\pc\geq 0$\\
    When Ann is observed to have $X=0$: \>$\pc\geq 0.75$\\
    When Ann's value for $X$ is not observed: \> $\pc\geq 0.71$\\[2ex]
    Ignoring $X$ and using the \textcite{tian/pearl:probcaus} bounds :
    \>$\pc\geq 0.53$.
  \end{tabbing}
\end{ex}

\section{Complete mediator}
\label{sec:compmed}
\begin{figure}[htbp]
  \centering
  \includegraphics[width=.5\linewidth]
  {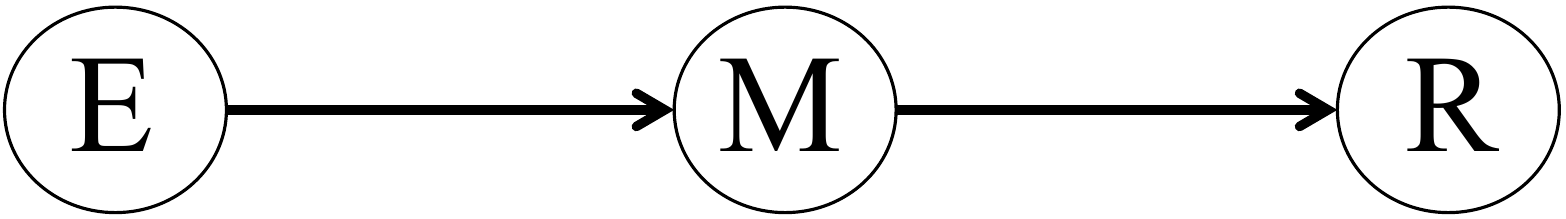}
  \caption{Complete mediator}
  \label{fig:compmed}
\end{figure}

We now turn to consider the case, represented by the diagram in
\figref{compmed}, that the causal effect of $E$ on $R$ is completely
mediated by some variable $M$, assumed for simplicity to be binary.
That is to say, $E$ affects $M$ and $M$ effects $R$, and the
combination of these two pathways constitutes the totality of the
effect of $E$ on $R$.
We assume that there is no further confounding of either of these
relationships, either in the data or for Ann.  We shall consider the
case that all three variables $(E,M,R)$ are observed in the data, but
that we only observe $E$ and $R$ for Ann.

Note that our assumptions imply the conditional independence property
$\ind R E M$, so that
\begin{equation}
  \label{eq:ab0}
  \pr(R=r \mid E = e) =  \sum_m\pr(R=r \mid M =m)\pr(M=m \mid E =e),
\end{equation}
which should hold, at least approximately, in the data.  In any event,
we should estimate \mbox{$\pr(R=r \mid E = e)$} using formula
\eqref{eq:ab0}, rather than directly.  (We could also have separate
experimental data on the $E\rightarrow M$ and $M\rightarrow R$
relationships, and apply \eqref{eq:ab0}.)

It is shown in \textcite{apd/rm/mm} that (with this understanding) the
additional information on the mediator does not affect the lower bound
in \eqref{eq:rrA}.  However, we obtain the following improved upper
bound:
\begin{eqnarray}
  \nonumber
  \lefteqn{\frac 1 {\pr(R=1 \mid E=1)}}\\[1.5ex]
  \nonumber
&{}\times\min \{&\pr(M=0 \mid E=0)\, \pr(R=0 \mid M=0)+\pr(R=0 \mid M=1)\,\pr(M=0 \mid E=1),\\
  \nonumber
&&{} \pr(M=1 \mid E=1)\, \pr(R=0 \mid M=0)+\pr(R=0 \mid M=1)\,\pr(M=1 \mid E=0),\\
  \nonumber
&&{}\pr(M=0 \mid E=0)\,\pr(R=1 \mid M=1)+\pr(R=1 \mid M=0)\,\pr(M=0 \mid E=1),\\
  \nonumber
&&{} \pr(M=1 \mid E=1)\,\pr(R=1 \mid M=1)+\pr(R=1 \mid M=0)\,\pr(M=1 \mid E=0)\quad\}.\\\quad
  \label{eq:compmedub}
\end{eqnarray}

\begin{ex}
  \label{ex:compmed}
  Suppose we obtain the following values from the data:
  \begin{eqnarray*}
    \Pr(M=1 \mid E= 1) &=& 0.25  \\
    \Pr(M=1 \mid E= 0) &=& 0.025 \\
    \Pr(R=1 \mid M= 1) &=& 0.9   \\
    \Pr(R=1 \mid M= 0) &=& 0.1.
  \end{eqnarray*}
  Marginalising over $M$ using \eqref{eq:ab0}, we obtain the same
  values as in \tabref{aspirin}.  On applying \eqref{eq:compmedub}, we
  get $0.60 \leq\pc \leq 0.76$; whereas without taking account of the
  mediator $M$ the bounds were $0.6 \leq \pc \leq 1$.
\end{ex}

\begin{figure}[htbp]
  \centering
  \includegraphics[width=.5\linewidth]
  {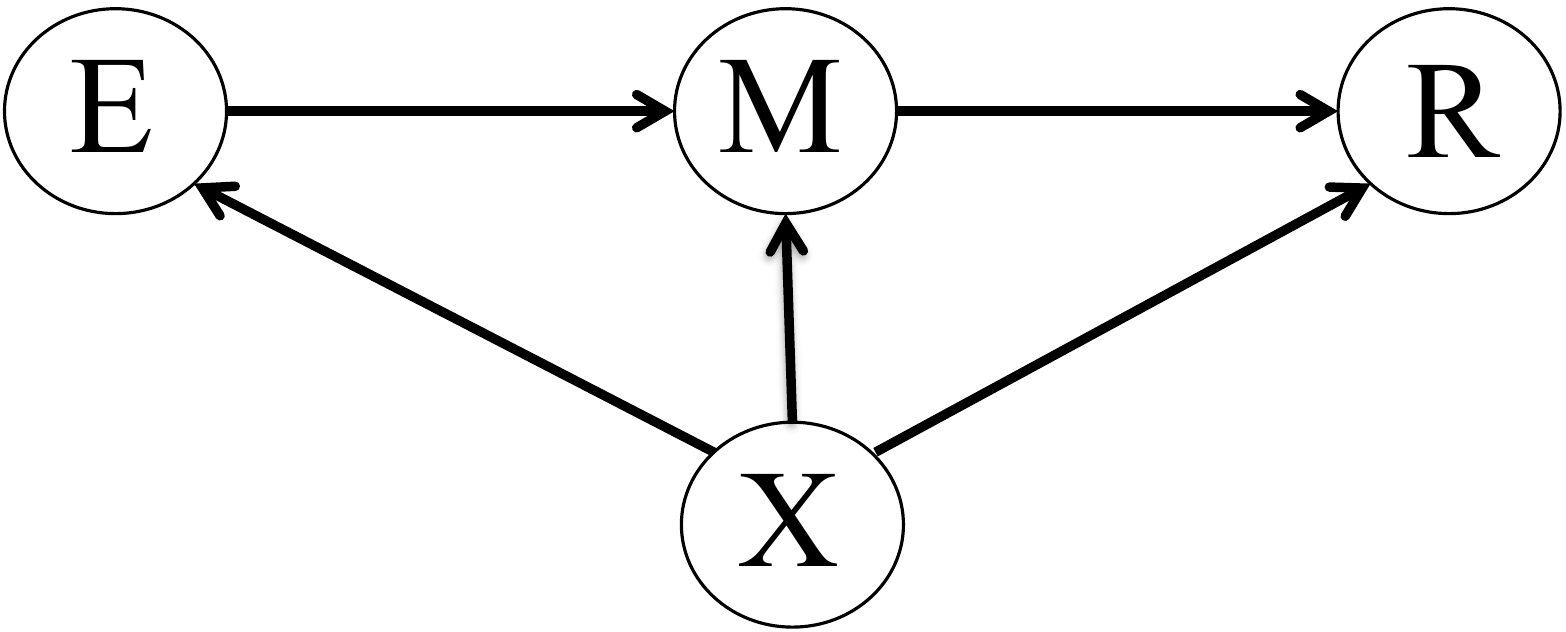}
  \caption{Complete mediator with covariate}
  \label{fig:allcomplete}
\end{figure}

We can elaborate the above analysis by allowing for an additional
sufficient covariate $X$, as represented by \figref{allcomplete}, that
can modify all the above marginal and conditional probabilities, and
conditional on which $M$ is a complete mediator.

We assume that $X$ is always observed in the data.  If we observe
$X=x$ for Ann, we simply replace all probabilities by their versions
further conditioned on $X=x$, yielding bounds
$l(x) \leq \pc \leq u(x)$, say.  If we do not observe $X$ in Ann, the
lower bound is $E\{l(X) \mid E=1, R=1)\}$, and the upper bound is
$E\{u(X) \mid E=1, R=1)\}$.

\section{Partial mediator}
\label{sec:partmed}

\begin{figure}[htbp] \centering
  \includegraphics[width=.5\linewidth]
  {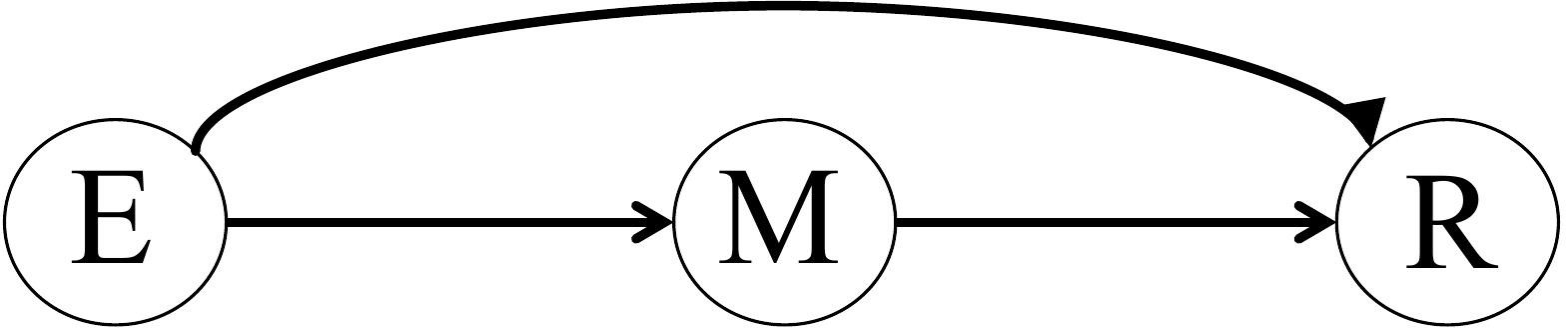}
  \caption{Partial mediator}
  \label{fig:partmed}
\end{figure}

In the case of a partial mediator, as represented by \figref{partmed},
there us an additional ``direct effect'' of $E$ on $R$, that is not
mediated by $M$ (again supposed binary, for simplicity).  So now we
have to consider the effect of $E$ on $M$, and the joint effect of $E$
and $M$ on $R$.  Again we suppose that there is no further confounding
of these relationships, either in the data or for Ann, and consider
the case that all three variables $(E,M,R)$ are observed in the data,
but that we only observe $E$ and $R$ for Ann.

In \textcite{rm/apd/mm} it is shown that, once again, the lower bound
for $\pc$ is unaffected by the additional information on $M$.
However, we are able to abtain a new upper bound $u$, where
\begin{eqnarray}
  \nonumber
  \pr(R=1 \mid E=1)\times u &=& 
                                \min\left\{\pr(R=0 \mid E= 0, M= 0),
                                \,\pr(R=1 \mid E= 1, M= 0)\right\}\\
  \nonumber
                            &&{}\,\,\,\,\times \min\left\{\pr(M=0 \mid E= 0),
                               \,\pr(M=0 \mid E= 1)\right\}\\
  \nonumber
                            &&+
                               \min\left\{\pr(R=0 \mid E= 0, M= 0),
                               \,\pr(R=1 \mid E= 1, M= 1)\right\}\\
  \nonumber
                            &&{}\,\,\,\,\times \min\left\{\pr(M=0 \mid E= 0),
                               \,\pr(M=1 \mid E= 1)\right\}\\
  \nonumber
                            &&+
                               \min\left\{\pr(R=0 \mid E= 0, M= 1),
                               \,\pr(R=1 \mid E= 1, M= 0)\right\}\\
  \nonumber
                            &&{}\,\,\,\,\times \min\left\{\pr(M=1 \mid E= 0),
                               \,\pr(M=0 \mid E= 1)\right\}\\      
  \nonumber
                            &&+
                               \min\left\{\pr(R=0 \mid E= 0, M= 1),
                               \,\pr(R=1 \mid E= 1, M= 1)\right\}\\
  \label{eq:partub}
                            &&{}\,\,\,\,\times \min\left\{\pr(M=1 \mid E= 0),
                               \,\pr(M=1 \mid E= 1)\right\}.
\end{eqnarray}

\begin{ex}
  \label{ex:partmed}
  Suppose that, from the data, we estimate the following
  probabilities:
  \begin{eqnarray*}
    \pr(M=1\mid E=0) &=& 0.27   \\
    \pr(M=1\mid E = 1) &=& 0.019 \\
    \pr(R=1\mid E=0, M=0) &=& 0.02  \\
    \pr(R=1\mid E=0, M=1) &=& 0.835  \\
    \pr(R=1\mid E=1, M=0) &=& 0.685 \\
    \pr(R=1\mid E=1, M=1) &=& 0.857.
  \end{eqnarray*}
  Marginalising over $M$, these imply
  \begin{eqnarray*}
    \pr(R=1\mid E=0) &=& 0.24  \\
    \pr(R=1\mid E=1) &=& 0.69.
  \end{eqnarray*}
  We then obtain $0.65 \leq \mbox{PC}\leq 0.81$ when accounting for
  the mediator, compared with $0.65 \leq \mbox{PC}\leq 1$ when
  ignoring it.
\end{ex}

However, accounting for a partial mediator may not always yield
improved bounds.
\begin{ex}
  Suppose the probabilities estimated from the data are as follows:
  \begin{eqnarray*}
    \pr(M=1 \mid E= 0) &=& 0.96\\
    \pr(M=1 \mid E= 1) &=& 0.74 \\
    \pr(R=1 \mid E= 0, M= 0)&=& 0.02  \\
    \pr(R=1 \mid E= 0, M= 1) &=& 0.33   \\
    \pr(R=1 \mid E= 1, M= 0) &=& 0.91   \\
    \pr(R=1 \mid E= 1, M= 1) &=& 0.73.
  \end{eqnarray*}
  These imply
  \begin{eqnarray*}
    \pr(R=1\mid E=0) &=& 0.32  \\
    \pr(R=1\mid E=1) &=& 0.78.
  \end{eqnarray*}
  The lower bound is $\pc \geq 0.59$.  The upper bound obtained from
  \eqref{eq:partub} is $\pc \leq 0.95$.  However, just ignoring $M$,
  and applying \eqref{eq:rrA}, we get upper bound $\pc \leq 0.88$,
  which should theregore be used in preference.
\end{ex}

\begin{figure}[htbp]
  \centering
  \includegraphics[width=.5\linewidth]
  {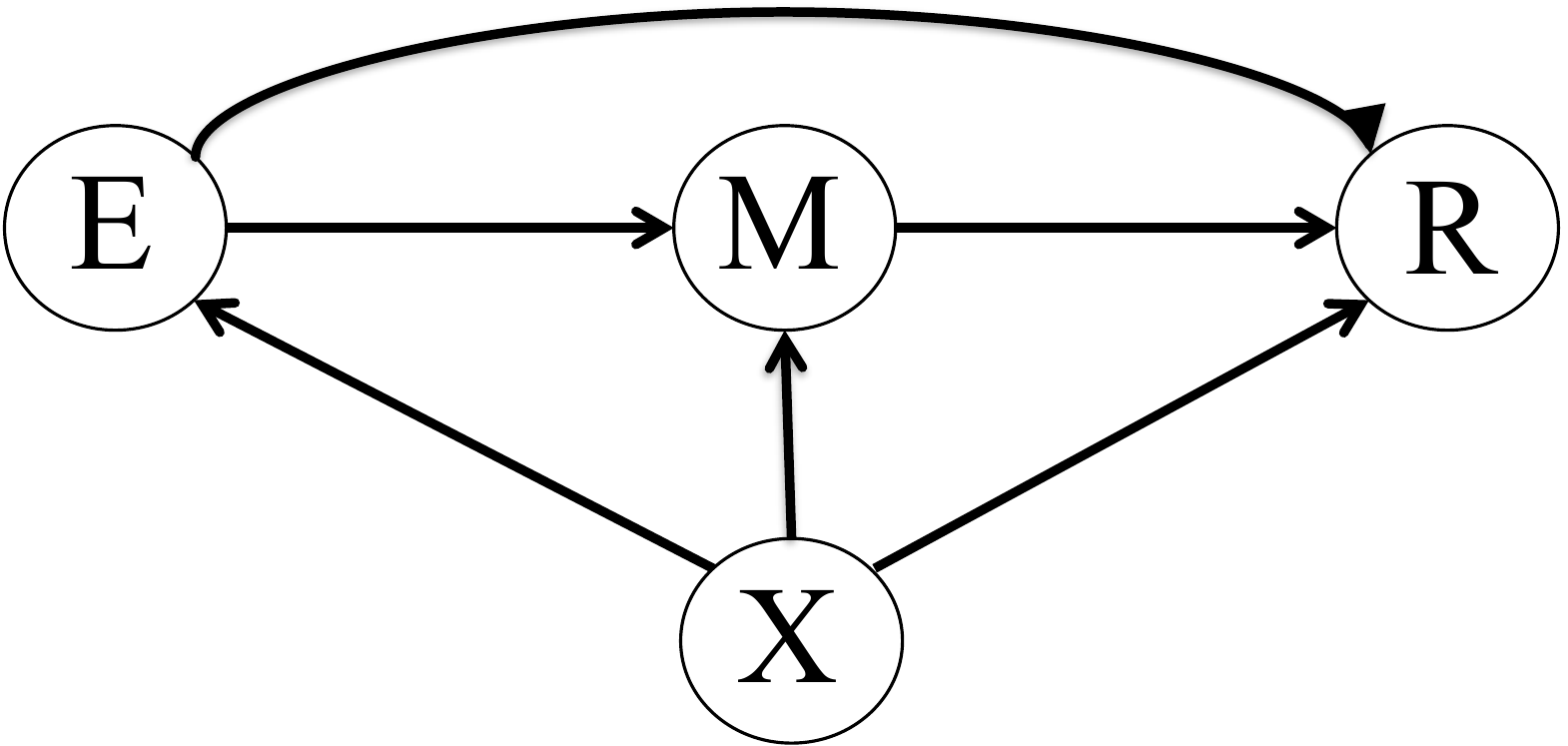}
  \caption{Partial mediator with covariate}
  \label{fig:allpartial}
\end{figure}

Again, we can elaborate the problem by including a sufficient
covariate $X$, conditional on which we have the above structure, as
represented by \figref{allpartial}.  Let the appropriate lower and
upper bounds, computed with all probabilities conditioned on $X=x$, be
$l(x)$ and $u(x)$.  Then in the absence of information about $X$ for
Ann, the relevant lower and upper bounds for $\pc$ are
$\E\{l(X) \mid E=1, R=1\}$, $\E\{u(X) \mid E=1, R=1\}$.

\section{Discussion}
\label{sec:disc}
We have assumed throughout that the available data are sufficiently
extensive to allow essentially perfect estimates of the required
probabilities.  Of course this will never be the case.
\textcite{apd/mm/sef:ba} consider how to make suitable inferences
about $\pc$ from limited data, using a Bayesian approach.

We have presented some simple cases, involving covariates and
mediators, where having additional information about the internal
structure of the causal link between exposure and response enables us
to refine our causal inferences.  Other problems might involve a more
complex collection of variables and pathways, perhaps expressed as a
directed acyclic graph.  It should be possible to extend our analysis
to such more general cases.  Information on different causal
subsystems might be obtained from a variety of scientific
investigations.  But caution must be exercised in justifying the
transfer of such information from generic scientific data to the
individual case.  Additional uncertainty over the correct
representation of the problem must also be taken into account.

To return, finally, to the issue facing a court of law that has to
return a decision on causality: we hope we have shown that this is not
a simple matter!  We can only express our sympathy with a judge who
has to make such a decision, on the ``balance of probabilities,'' when
the best possible conclusion, based on available scientific evidence
and the kind of analysis we have presented here, is that the
probability of causation lies between $0.4$ and $0.8$.




\end{document}